\documentclass[12pt]{amsart}
\usepackage{amssymb}

\begin{document}
\newcommand{\nt}{\noindent}
\newcommand{\bs}{\medskip}
\newcommand{\ms}{\smallskip}
\newcommand{\mk}{\smallskip}
\newcommand{\sk}{\smallskip}
\newcommand{\ve}{\varepsilon}
\newcommand{\vf}{\varphi}

\title{Vanishing Sums of $\,m\,$th Roots of Unity in Finite Fields}

\author{T. Y. Lam}
\address{\hskip-\parindent
T. Y. Lam, Mathematics Department, University of California,
Berkeley, CA, 94720}
\email{lam@msri.org}

\thanks{Lam was supported in part by NSF.
Research at MSRI is supported in part by NSF grant DMS-9022140.}

\author{K.~H.~Leung}
\thanks{Leung's research was carried out while he
was on sabbatical leave at U.C.~Berkeley from the National University
of Singapore. The hospitality of the former institution is gratefully
acknowledged.}

\address{\hskip-\parindent K. H. Leung, National University of
Singapore, Singapore 119260}

\begin{abstract}
\begin{small}
In an earlier work, the authors have determined all possible weights
$\,n\,$ for which there exists a vanishing sum $\,\zeta_1+\cdots +\zeta_n=0\,$
of $\,m\,$th roots of unity $\,\zeta_i\,$ in characteristic 0.  In this paper,
the same problem is studied in finite fields of characteristic $\,p$.  
For given $\,m\,$ and $\,p$, results are obtained on integers $\,n_0\,$ 
such that all integers $\,n\geq n_0\,$ are in the ``weight set'' $\,W_p(m)$. 
The main result $\,(1.3)\,$ in this paper guarantees, under suitable 
conditions, the existence of solutions of $\,x_1^d+\cdots+x_n^d=0\,$
{\it with all coordinates not equal to zero\/} over a finite field.
\end{small}
\end{abstract}

\maketitle

\section{Introduction}

By a vanishing sum of $\,m\,$th roots of unity, we mean an equation
$\,\alpha_1+\cdots+\alpha_n=0\,$ where $\,\alpha_i^m=1\,$ for each $\,i$.
The integer $\,n\,$ is said to be the {\it weight\/} of this vanishing
sum.  In [LL], considering $\,m\,$th roots of unity in $\,{\mathbb C}$,
we defined $\,W(m)\,$ to be the set of integers $\,n\geq 0\,$ for which there 
exists a vanishing sum $\,\alpha_1+\cdots+\alpha_n=0\,$ as above.  The
principal result in [LL] gives a complete determination of the weight set
$\,W(m)\,$ (in characteristic 0), as follows.

\bs
\nt {\bf Theorem 1.1.} {\it For any natural number $\,m\,$ with prime
factorization $\,p_1^{a_1}\cdots p_r^{a_r}$, the weight set $\,W(m)\,$ is
exactly given by} $\,{\mathbb N}\,p_1+\cdots+{\mathbb N}\,p_r.$ ({\it Here and in
the following,} $\,{\mathbb N}:=\{0,1,2,\cdots\})$.

\bs
In this paper, we study vanishing sums of $\,m\,$th roots of unity {\it in
characteristic $\,p$.}  In analogy to the characteristic $\,0\,$ case, we 
define $\,W_p(m)\,$ to be the set of weights $\,n\,$ for which there exists a
vanishing sum $\,\alpha_1+\cdots+\alpha_n=0\,$ where each $\,\alpha_i\,$ is an
$\,m\,$th root of unity in $\,\overline{{\mathbb F}}_p$, the algebraic closure
of the prime field $\,{\mathbb F}_p$. Note that, if $\,m=p^tm'\,$ where
gcd$(p,\,m')=1$, we have $\,x^m=1\,$ in $\,\overline{{\mathbb F}}_p\,$ iff
$\,x^{m'}=1$; in particular, $\,W_p(m)=W_p(m')$.  Therefore, we may assume
throughout that $\,\mbox{gcd}(p,\,m)=1$, i.e. $\,p\,$ is not among the 
prime divisors $\,p_i\,$ of $\,m$.
As in the case of charactersitic 0, we have $\,p_i\in W_p(m)\,$ for all $\,i$.
But in characteristic $\,p$, we also have $\,p\in W_p(m)\,$ (due to the
vanishing sum $\,p\cdot 1=0$), so now
$$ W_p(m)\supseteq {\mathbb N}\,p+{\mathbb N}\,p_1+\cdots+{\mathbb N}\,p_r.
\leqno (1.2)$$
\indent Easy examples  (see (2.1)) show that this need not be an equality in
general, so we are left with no viable conjecture on the structure of the 
weight set $\,W_p(m)\,$ in characteristic $\,p$.  However, (1.2) does show 
that, if $\,m>1$, all sufficiently large integers $\,n\,$ (in fact all
$\,n\geq (p-1)(p_i-1)$) belong to $\,W_p(m)$.  A more tractable problem 
will then be the determination of more accurate bounds $\,n_0\,$ such that 
all integers $\,n\geq n_0\,$ belong to $\,W_p(m)$. 

\bs
In this paper, we will show how such an integer $\,n_0\,$ can be determined. 
Our work is divided into three cases, depending on whether
$\,\mbox{gcd}(p-1,\,m)\,$ is $\,1,\,2$, or bigger.  The estimates
on $\,n_0\,$  differ from case to case, and are given respectively in (5.6), 
(4.1), and (3.1)-(3.3).  Although we have three different estimates on
$\,n_0$, there does exist a (necessarily weaker) uniform estimate for all 
cases.  In the following, we shall try to explain what this uniform estimate 
is, and why is it a reasonable one.

\bs
A guiding principle for our work throughout is the fact that a finite field is
a $\,C_1$-field (see [Gr]).  If $\,K={\mathbb F}_{p^k}\,$ is a finite field 
containing all $\,m\,$th roots of unity, then, for $\,d:=(p^k-1)/m$, the 
$\,m\,$th roots of unity in $\,\overline{\mathbb F}_p\,$ comprise the group 
$\,\dot{K}^{\,d}$.  Therefore, a vanishing sum of $\,m\,$th roots
 of unity of weight $\,n\,$ corresponds precisely to a ``good'' solution of 
$\,x_1^d+\cdots+x_n^d=0\,$ in $\,K$, where by a ``good'' solution
we mean one with {\it each\/} $\,x_i\neq 0$.  If
$\,n>d$, the fact that $\,K\,$ is $\,C_1\,$ implies that
we have a solution $\,(x_1,\,\cdots,\,x_n)\neq (0,\,\cdots,\,0)$.  It certainly
seems tempting to speculate that there exists in fact a 
``good'' solution (in $\,K$).  If this is indeed the case, then 
by what we said earlier in this paragraph, {\it any\/} integer $\,n>d\,$ will 
be in the weight set $\,W_p(m)$.

\bs
The desired conclusion that, for $\,n>d$, $\,x_1^d+\cdots+x_n^d=0\,$ has a
``good'' solution in $\,K\,$ is, however, not true in general!
For instance, if $\,d=p^k-1$, then $\,x^d=1\,$ for each $\,x\in \dot{K}$,
so we have a ``good'' solution for $\,x_1^d+\cdots+x_n^d=0\,$ in $\,K\,$
only when $\,n\,$ is a multiple of $\,p$.  In a similar vein,
if $\,k=1$, $\,p\,$ is odd, and $\,d=(p-1)/2$, then any nonzero 
$\,d\,$th power in $\,K\,$ is $\,\pm 1$.  For any odd integer $\,n\in (d,\,p)$,
the equation $\,x_1^d+\cdots+x_n^d=0\,$ again has no ``good'' solution in 
$\,{\mathbb F}_{p}$.  The trouble with these cases is that $\,m\leq 2$,
for which we don't have ``enough'' $\,m\,$th roots to play with.  As it
turns out, as soon as we ignore the above cases, we'll have the 
following uniform result for getting ``good'' solutions.

\bs
\nt {\bf Theorem 1.3.} {\it Let} $\,K={\mathbb F}_{p^k}\,$ {\it and 
$\,d=(p^k-1)/m\,$ as above, and assume that $\,m\neq 1$, $\,(m,\,k)\neq (2,1)$.
Then, whenever $\,n>d$, the equation
$\,x_1^d+\cdots+x_n^d=0\,$ has a ``good'' solution in $\,K$.  In
other words, the weight set $\,W_p(m)\,$ contains all integers $\geq d+1$. }

\bs
The results in \S\S3-5 below will cover this theorem in the case $\,m\geq 3$.
In the case $\,m=2$, (1.3) is quickly checked as follows.  Since we assume 
in this case that $\,k\geq 2$, we have $\,d\geq (p^2-1)/2\geq p-1$.  Given
$\,n\geq d+1\geq p$, it is easy to solve the equation 
$\,\alpha_1+\cdots+\alpha_n=0\,$ with $\,\alpha_i=\pm 1$, by considering
the parity of $\,n$.  Having disposed of the trivial cases $\,m=1,\,2$, 
{\it we may assume in\/} \S\S3-6 {\it of this paper that $\,m\geq 3$.}

\bs
In the case when $\,p\,$ is odd and $\,d=2$, (1.3) says precisely that, for
any $\,n>2$, the quadratic form $\,x_1^2+\cdots+x_n^2\,$ has a ``good'' zero 
over any finite field of more than five elements.  This is a special case of
a well-known observation of Witt for isotropic diagonal quadratic forms (see 
[Wi:~p.39], or [BS:~p.394], [La:~p.25, Ex.7]).  Thus, (1.3) may be thought of 
as a generalization of Witt's result to the higher degree diagonal forms
$\,x_1^d+\cdots+x_n^d\,$ over finite fields.  Note that $\,d\,|\,(p^k-1)\,$
is not a really essential assumption in (1.3).  In dealing with the equation
$\,x_1^e+\cdots+x_n^e=0$, we can replace the degree $\,e\,$ by
$\,d:=\mbox{gcd}(p^k-1,\,e)$, and define $\,m\,$ to be $\,(p^k-1)/d$.
Then $\,\dot{K}^e=\dot{K}^d$, so under the assumptions of (1.3),
$\,x_1^e+\cdots+x_n^e=0\,$ will have a ``good'' solution as long as $\,n>d$.

\bs
In the literature, there are many results dealing with diagonal equations
over finite fields; see, for instance, [LN], [Sch], [Sm], and more recently
[QY].  Conventionally, one could apply algebro-geometric methods, or
alternatively the method of Gauss and Jacobi sums.  As the referee of this
paper pointed out, these methods can be utilized to show the existence of 
``good'' solutions to a diagonal equation $\,x_1^d+\cdots+x_n^d=0\,\;(n>2)\,$ 
in $\,{\mathbb F}_q\,$ if $\,q\,$ is suitably large compared to $\,d\,$ 
(without the condition $\,n>d\,$).  However, these conventional methods do
not seem to give enough information if $\,q\,$ is ``small'' in comparison
to $\,d$.  In our setting, working mostly with $\,n>d\,$ and taking full
advantage of the additive nature of the special equation 
$\,x_1^d+\cdots+x_n^d=0$, we apply instead the methods of additive 
number theory.  These methods do give fairly precise results, without
reference to the size of $\,{\mathbb F}_q\,$.  In fact, the analysis in 
\S\S3-5 will not only prove (1.3), but also show that, in various cases, 
the equation $\,x_1^d+\cdots+x_n^d=0\,$ has a ``good'' solution in 
$\,K\,$ often for much smaller values of $\,n\,$ (than $\,n\geq d+1$).  
Thus, the more precise results in this paper are to be found in (3.1)-(3.3), 
(4.1) and (5.6).  Theorem 1.3 is only a common denominator of these results 
giving a convenient and uniform summary of the main work in this paper.

\bs
\noindent {\it Acknowledgment.}  We thank the referee of this
paper, whose comments enabled us to rewrite more accurately the last 
paragraph above comparing the use of different methods in treating
diagonal equations over finite fields.

\section{Some Basic Examples}

We shall begin with some examples and computations of the weight sets
$\,W_p(m)$.  The first couple of examples show that various properties of 
weight sets in characterictic 0 are no longer valid in characteristic $\,p$.
For convenience of expressing weight sets, let us use the notation 
$\,[\,n,\,\infty)_{{\mathbb Z}}\,$ for the set of integers $\geq n$. 

\bs
\nt {\bf Example 2.1.} Referring to (1.2), {\it the smallest positive element 
in the set $\,W_p(m)\,$ may not be\/}
$\,\mbox{min}\,\{p,\,p_1,\cdots,p_r\}$.  For instance,
when $\,p=11\,$ and $\,m=5$, the $\,5$th roots of unity in 
$\,\overline{{\mathbb F}}_{11}\,$ are $\,\{1,3,9,5,4\}$.  Observing that
$\,1+1+9=0\,$ in $\,{\mathbb F}_{11}$, we see that $\,W_{11}(5)\,$ contains
$\,3$, which is smaller than $\,5\,$ and $\,11$.  By (2.3) below, we have
$\,W_{11}(5)=\{0\}\cup [\,3,\,\infty)_{{\mathbb Z}}$.  Thus, not only (1.2) fails
to be an equality, but also $\,W_{11}(5)\,$  is not 
even of the form $\,\sum_i\,{\mathbb N}\,q_i\,$ for a set of primes $\,q_i\,$'s.

\bs
\nt {\bf Example 2.2.}  Contrary to the characteristic $0$ case,
{\it the set $\,W_p(m)\,$ may be larger than $\,W_p(m_0)\,$ where 
$\,m_0\,$ is the square-free part of $\,m$.}  For instance, let $\,p=5\,$ 
and $\,m=4$, so $\,m_0=2$.  It is easy to see that 
$\,W_5(2)=\{0,\,2\}\!\cup\![\,4,\,\infty)_{{\mathbb Z}}$, but
$\,W_5(4)=\{0\}\cup [\,2,\,\infty)_{{\mathbb Z}}$.

\bs
\nt {\bf Example 2.3.} Let $\,q=p^a>5\,$ where $\,p\,$ is an odd prime,
and let $\,m=(q-1)/2$. Then $\,d:=(q-1)/m=2$. For any $\,n\geq 3$, the
quadratic form $\,X_1^2+\cdots+X_n^2\,$ is isotropic over $\,{\mathbb F}_q$,
so by the theorem of Witt referenced before, it has a ``good'' zero in 
$\,{\mathbb F}_q$.  Therefore, $\,n\in W_p(m)$.  It follows that
$\,W_p(m)=\{0\}\cup [\,2,\,\infty)_{{\mathbb Z}}\,$ 
     if $\,q\equiv 1\,(\mbox{mod}\;4)$,  and
$\,W_p(m)=\{0\}\cup [\,3,\,\infty)_{{\mathbb Z}}\,$ 
     if $\,q\equiv 3\,(\mbox{mod}\;4)$. 

\bs
\nt {\bf Example 2.4.} (${\mathbb F}_p\,$ contains all $\,m\,$th roots of unity.)
Let $\,p=31$, and $\,m=3$.  The third roots of unity are $\,\{1,5,25\}$,
so the equation $\,25+6\cdot 1=0\in {\mathbb F}_{31}\,$ shows that
$\,7\in W_{31}(3)$.  A routine computation shows that
$\,W_{31}(3)=\{0,\,3,\,6,\,7\}\cup [\,9,\,\infty)_{{\mathbb Z}}$.

\bs
\nt {\bf Example 2.5.} (${\mathbb F}_p\,$ contains no $\,m\,$th roots of unity
other than 1.)
Let $\,p=2$, and $\,m=73$. We work in $\,K={\mathbb F}_{2^9}\,$ which contains
all $\,73\,$rd roots of unity.  By standard tables of irreducible polynomials 
over finite fields, the trinomial
$\,f(X)=X^9+X+1\,$ is irreducible over $\,{\mathbb F}_2$, so we can take
$\,K\,$ to be $\,{\mathbb F}_2[X]/(f(X))$. Let $\,\alpha:=\overline{X}\in K$.
We have $\,0=(\alpha^9+\alpha+1)^8=\alpha^{72}+\alpha^8+1$,
so $\,\alpha^{73}=\alpha(\alpha^8+1)=\alpha^9+\alpha=1$. Thus,
the relation $\,\alpha^9+\alpha+1=0\,$ shows that $\,3\in W_2(73)$, 
and it follows easily that $\,W_2(73)=\{0\}\cup [\,2,\,\infty)_{{\mathbb Z}}$.

\bs
In the balance of this section, let us consider $\,W_p(m)\,$ in the case
when $\,m\,$ is a prime power (not divisible by $\,p$).  Under a special 
hypothesis on the cyclotomic polynomial $\,\Phi_m(X)$, the weight set 
$\,W_p(m)\,$ can be determined explicitly. 

\bs
\nt {\bf Theorem 2.6.} {\it Let $\,m=\ell\,^a\,$ where $\,\ell\,$ is a prime
different from $\,p$, and assume that the cyclotomic polynomial}
$\,\Phi_m(X)\in {\mathbb Z}[X]\,$ {\it remains irreducible modulo $\,p$.
Then} $\,W_p(m)={\mathbb N}\,p+{\mathbb N}\,\ell$.

\bs
\nt {\bf Proof.} Of course, it suffices to prove the inclusion ``$\subseteq$''.
Let $\,\zeta\,$ be a primitive $\,m\,$th root of unity
in $\,\overline{{\mathbb F}}_p$. Let $\,m':=m/\ell$, and $\,\alpha:=\zeta^{m'}$
(a primitive $\,\ell\,$th root of unity).  Let $\,K={\mathbb F}_p(\zeta)$,
and $\,L={\mathbb F}_p(\alpha)$.  Since $\,\Phi_m(X)\,$ is irreducible 
$\,\mbox{mod}\;p$, $\,[K:{\mathbb F}_p]=\varphi(m)=m'(\ell-1)$.  From this,
it is easy to see that $\,[K:L]=m'\,$ and $\,[L:{\mathbb F}_p]=\ell-1$.

\ms
Any vanishing sum of $\,m\,$th roots of unity can be written in the form
$\,\sum_{i=0}^{m'-1}\,g_i\zeta^i=0$, where each $\,g_i\,$ is a sum of
$\,\ell\,$th roots of unity.  Since the degree of $\,\zeta\,$ over $\,L\,$
is $\,m'$, the elements $\,1,\,\zeta,\,\cdots\,,\,\zeta^{m'-1}\,$
are linearly independent over $\,L$.  Therefore, each $\,g_i\in L\,$ is 
itself a vanishing sum, and it suffices to show that its weight is in
$\,{\mathbb N}\,p+{\mathbb N}\,\ell$.  Starting over again, we are now down to
considering  a vanishing sum $\,\sum_{i=0}^{\ell-1}\,a_i\alpha^i=0$, where
each $\,a_i\in {\mathbb N}$.  Let $\,a_j\,$ be the smallest among the $\,a_i\,$'s.
Since the minimal equation of $\,\alpha\,$ over $\,{\mathbb F}_p\,$ is
$\,1+\alpha+\cdots+\alpha^{\ell-1}=0$, it follows easily that
$\,a_0=a_1=\cdots=a_{\ell-1}\in {\mathbb F}_p$. The weight of the vanishing sum 
in question is $\,\sum_i\,a_i\equiv \ell\,a_j \;(\mbox{mod}\;p)$.
Since $\,\sum_i\,a_i\geq \ell\,a_j$, it follows that
$\,\sum_i\,a_i=\ell\,a_j+bp\,$ for some $\,b\in {\mathbb N}$, as desired. 
\qed

\bs
\nt {\bf Remark 2.7.} A vanishing sum of $\,m\,$th roots of unity is said
to be {\it minimal\/} if no proper subsum of it is also vanishing. In general,
the problem of determining the minimal vanishing sums is difficult (both in
characteristic 0 and in characteristic $\,p$). Under the hypothesis of (2.6), 
however, this problem can be solved.  In fact, the argument presented in the
proof above can be used to show that, in the setting of (2.6), the minimal
vanishing sums of $\,m\,$th roots of unity are, up to multiplication by
a power of $\,\zeta$: 
(1) $\,p\cdot 1=0$, and (2) $\,1+\alpha+\cdots+\alpha^{\ell-1}=0$.  
(Of course, this implies that $\,W_p(m)={\mathbb N}\,p+{\mathbb N}\,\ell$.)
For this conclusion, however, the 
assumption on the irreducibility of $\,\Phi_m(X)\,$ modulo $\,p\,$ is 
essential, as the examples (2.1), (2.4) and (2.5) show.  (In (2.1) and
(2.4), $\,\Phi_m(X)\,$ splits completely modulo $\,p$, and in (2.5),
$\,\Phi_m(X)\,$ splits into the product of eight irreducible
factors of degree 9 in $\,{\mathbb F}_p[X]$.)

\section{The Case $\,\mbox{gcd}(p-1,\,m)\geq 3$}

In dealing with $\,W_p(m)$, our main goal is to find good estimates for 
integers $\,n_0\,$ such that $\,[n_0,\,\infty)_{{\mathbb Z}}\subseteq W_p(m)$.
We begin our analysis with the case
when $\,\mbox{gcd}(p-1,\,m)\geq 3$.  This case turns out to be fairly easy if
we use the right tools from additive number theory modulo $\,p$.  It will
be convenient to use the following notations.  For a subset $\,A\,$ in a field,
we shall write $\,|A|\,$ for the cardinality of $\,A$, and for any integer
$\,n\geq 1$, we write $\,n\ast A\,$ for the set $\,A+\cdots+A\,$ with
$\,n\,$ summands of $\,A$.

\bs
\nt {\bf Theorem 3.1.} {\it Assume that}
$\,m_0:=\mbox{gcd}(p-1,\,m)\geq 3\,$ {\it and let $\,d_0=(p-1)/m_0$.
Then $\,[\,d_0+1,\,\infty)_{{\mathbb Z}}\subseteq W_p(m_0)\subseteq W_p(m)$.}

\bs
\nt {\bf Proof.} Since $\,m_0|(p-1)$, the group $\,H\,$ of $\,m_0\,$th roots
of unity in $\,{\mathbb F}_p\,$ is exactly $\,\dot{{\mathbb F}}^{\,d_0}_p\,$ and has
exactly $\,m_0\,$ elements.  We claim that
$\,|n\ast H|\geq nm_0\,$ for $\,n\leq d_0$, and
$\,|n\ast H|=p\,$ for $\,n\geq d_0+1$.  It suffices to prove this for 
$\,n=1,2,\,\cdots,\,d_0+1\,$ (for, once we show that
$\,|(d_0+1)\ast H|=p$, then $\,(d_0+1)\ast H={\mathbb F}_p$, and this implies that
$\,(d_0+i)\ast H={\mathbb F}_p\,$ for any $\,i\geq 1$).  We proceed by induction
on $\,n$, the case $\,n=1\,$ being clear.  Assume that
$\,|n\ast H|\geq nm_0\,$ where $\,n<d_0$.  By the Cauchy-Davenport Theorem
(see [Ma: Cor.1.2.3]), $\,|(n+1)\ast H|\,$ is either $\,p\,$ (and hence
$\geq (n+1)m_0$), or else
$$ |(n+1)\ast H|\geq |n\ast H|+|H|-1\geq (n+1)m_0-1. $$
In the latter case, $\,|(n+1)\ast H\setminus\{0\}|\geq (n+1)m_0-2$.
Since $\,H\,$ acts on $\,(n+1)\ast H\setminus \{0\}\,$ by multiplication,
$\,|(n+1)\ast H\setminus \{0\}|\,$  is a multiple of $\,m_0$.  Since
$\,m_0\geq 3$, we must therefore have 
$\,|(n+1)\ast H\setminus \{0\}|\geq (n+1)m_0$,
which gives what we want.  This proves our claim for $\,n\leq d_0$.
In particular, $\,|d_0\ast H|\geq d_0m_0=p-1$.  By the Cauchy-Davenport
Theorem again, $\,(d_0+1)\ast H\,$ must be $\,{\mathbb F}_p$, for otherwise 
we would have
$$ |(d_0+1)\ast H|\geq |d_0\ast H|+|H|-1\geq d_0m_0+m_0-1=p+(m_0-1)>p, $$
a contradiction.  This completes our inductive proof.  Thus, for any
$\,n\geq d_0+1$, we have $\,0\in n\ast H$.  This means that
$\,n\in W_p(m_0)$, and so 
$\,[d_0+1,\,\infty)_{{\mathbb Z}}\subseteq W_p(m_0)\subseteq W_p(m)$.
\qed

\bs
\nt {\bf Example 3.2.} In many cases Theorem 3.1 gives the best result.
For instance, if $\,p\equiv 3\,(\mbox{mod}\;4)\,$ and $\,m=(p-1)/2\geq 3$,
then $\,d_0=2\,$ and we have $\,d_0\notin W_p(m)\,$ since $\,m\,$ is odd.
Even in the case $\,p\equiv 1\,(\mbox{mod}\;4)$, the Theorem may still give
the best result.  For instance, if $\,p=13\,$ and $\,m=4$, then $\,d_0=3\,$
and $\,G=\{\pm 1,\,\pm 8\}$.  By a simple calculation, 
$\,3\ast G=\dot{{\mathbb F}}_{13}$, so again $\,d_0=3\notin W_p(m)$.  On the other
hand, if $\,m\,$ is divisible by two distinct primes $\,p_1,\;p_2$, then the
fact that $\,p_1,\,p_2\in W_p(m)\,$ implies that
$\,[n_0,\,\infty)_{{\mathbb Z}}\subseteq W_p(m)\,$ for $\,n_0:=(p_1-1)(p_2-1)\,$ 
(see [LeV:~p.22, Ex.4]).  In case the number $\,d_0\,$ in (3.1) is ``large'',
$\,[n_0,\,\infty)_{{\mathbb Z}}\subseteq W_p(m)\,$ will of course give a better
result.

\bs
We can now derive the first case of Theorem 1.3.

\bs
\nt {\bf Corollary 3.3.} {\it Let $\,K={\mathbb F}_{p^k}\,$ ba a finite field
containing all $\,m\,$th roots of unity, and let $\,d=(p^k-1)/m$.  If}
$\,m_0:=\mbox{gcd}(p-1,\,m)\geq 3$, {\it then 
$\,[d+1,\,\infty)_{{\mathbb Z}}\subseteq W_p(m_0)\subseteq W_p(m)$.}

\bs
\nt {\bf Proof.} Say $\,p-1=m_0d_0\,$ and $\,m=m_0m_1$.  Then
$$ d=\frac{(p-1)(p^{k-1}+\cdots+p+1)}{m}=d_0\cdot \frac{p^{k-1}+\cdots+p+1}{m_1}.  $$
Since $\,\mbox{gcd}(d_0,\,m_1)=1$, the fraction on the RHS above is an integer.
Therefore, we have $\,d_0\,|\,d$, and the desired conclusion follows from
Theorem 3.1. \qed

\section{The Case $\,\mbox{gcd}(p-1,\,m)=2$}

We shall assume throughout this section that $\,\mbox{gcd}(p-1,\,m)=2\,$
(and as before $\,m\geq 3$).  In particular, $\,p\,$ is odd and 
$\,m\,$ is even.  In this case, $\,W_p(m)\,$
contains $\,2\,{\mathbb N}\,$ and is stable under addition by $\,2$.  Thus, once
we have an odd integer $\,n\in W_p(m)$, we will have automatically
$\,[n-1,\,\infty)_{{\mathbb Z}}\subseteq W_p(m)$.  This observation will be
used without further mention in the following.

\bs 
Let $\,K={\mathbb F}_{p^k}\,$ be any finite field containing the group
$\,G\,$ of all $\,m\,$th roots of unity.  The following result gives a
somewhat sharper form of Theorem 1.3 in the case 
$\,\mbox{gcd}(p-1,\,m)=2\,$ (in that the index $\,[\dot{K}:G]\,$ itself is
shown to be a weight, with a minor exception).

\bs
\nt {\bf Theorem 4.1.} {\it Assume that}
$\,\mbox{gcd}(p-1,\,m)=2$,
{\it and let $\,d=[\dot{K}:G]=(p^k-1)/m$.  Then
$\,[\,d,\,\infty)_{{\mathbb Z}}\subseteq W_p(m)\,$ unless $\,p=3\,$ and
$\,m=3^k-1\,$, in which case
$\,[\,d+1,\,\infty)_{{\mathbb Z}}\subseteq W_3(m)$. }

\bs
\nt {\bf Proof.} Let us first check the Theorem when $\,m=4$. In this case,
the assumption  $\,\mbox{gcd}(p-1,\,m)=2\,$ implies that $\,G\,$ is not
contained in $\,{\mathbb F}_p$, so $\,k\geq 2$.  If $\,p>3$, then
$\,d\geq (p^2-1)/4\geq p$, and we have 
$\,[\,d,\,\infty)_{{\mathbb Z}}\subseteq [\,p,\,\infty)_{{\mathbb Z}}\subseteq W_p(m)$.
If $\,p=3$, then $\,d\geq (9-1)/4=2$, and we have again
$\,[\,d,\,\infty)_{{\mathbb Z}}\subseteq [\,2,\,\infty)_{{\mathbb Z}}\subseteq W_3(m)$
(since $\,W_3(m)\,$ contains both $\,2\,$ and $\,3$).
{\it In the following, we may therefore assume that $\,m\geq 6$.}

\ms
Write $\,m=2m'$, so  that
$$ d=\frac{(p-1)(p^{k-1}+\cdots+p+1)}{2m'}=\frac{p-1}{2}\cdot \frac{p^{k-1}+\cdots+p+1}{m'}.  \leqno (4.2)  $$
Since $\,\mbox{gcd}(m',\,(p-1)/2)=1$, we have $\,m'|(p^{k-1}+\cdots+p+1)$.
If $\,m'<p^{k-1}+\cdots +p+1$, the second factor on the RHS in (4.2) is
$\,\geq 2$, so $\,d\geq p-1$.  Since $\,p\in W_p(m)$, we have
$\,[d,\,\infty)_{{\mathbb Z}}\subseteq [p-1,\,\infty)_{{\mathbb Z}}\subseteq W_p(m)$,
as desired.  Therefore, {\it in the following we may assume that}
$$\,m'=p^{k-1}+\cdots+p+1, \;\;\, \mbox{and} \,\;\;d=(p-1)/2.  \leqno (4.3)  $$
In this case $\,\dot{K}=\dot{{\mathbb F}}_{p-1}\!\cdot\!G$, so any coset of 
$\,G\,$ in $\,\dot{K}\,$ has a ``scalar'' representative. We fix a generator 
$\,\zeta\,$ for the group $\,G$, and try to put a lower bound on the 
cardinality of the set $\,A:={\mathbb F}_p\cap 2\ast G$. 

\bs
Recalling that $\,m\geq 6$,  write 
$$  (\zeta-1)G=a_1G\,,\;\;\,  (\zeta^2-1)G=a_2G\,, \;\;\,
\mbox{and} \;\;\; (\zeta^4-1)G=a_3G\,\;\;\,\mbox{where} \;\;
a_i\in \dot{{\mathbb F}}_p.  $$
Clearly, $\,\pm a_i\in A$, since
$\,-1\in G$.  {\it First let us assume that these three $\,G$-cosets in 
$\,\dot{K}\,$ are different.}  Since $\,a_i\neq -a_j$, $\,\{0,\,\pm a_i\}\,$
are seven different elements of $\,A$.  (In particular, $\,p\geq 7\,$ here.)
Applying repeatedly the Cauchy-Davenport Theorem in $\,{\mathbb F}_p$, we see that
$\,|n\ast A|\geq \mbox{min}\{p,\,6n+1\}$.  It follows that
$$ n\geq (p-1)/6 \Longrightarrow n\ast A={\mathbb F}_p
                 \Longrightarrow -1\in n\ast A
                 \Longrightarrow 2n+1\in W_p(m).  $$
This yields $\,2\lceil \frac{p-1}{6}\rceil+1 \in W_p(m)\,$ (where 
$\,\lceil \cdot \rceil\,$ denotes the ceiling function).  Writing $\,p\,$ 
in the form $\,6t\pm 1$, we see easily that $\,2\lceil \frac{p-1}{6}\rceil
=\lceil \frac{p-1}{3}\rceil $.  Thus, in this case, we get the stronger 
conclusion that $\,[\,(p-1)/3,\,\infty)_{{\mathbb Z}}\subseteq W_p(m)$.

\ms
{\it From now on, we may assume that the three cosets $\,\{a_iG\}\,$ above
are not all different}.  If $\,a_1G=a_2G$, then 
$\,\zeta^2-1=(\zeta-1)\zeta^i\,$ for some $\,i$, and so $\,\zeta^i=\zeta+1$.
Since $\,m\,$ is even, this shows that $\,3\in W_p(m)$, and so 
$\,[\,2,\,\infty)_{{\mathbb Z}}\subseteq W_p(m)$.  We have certainly no problem 
in this case (except when $\,d=1$, which occurs only when $\,p=3$).
If $\,a_2G=a_3G$, we can finish similarly.  Now assume $\,a_1G=a_3G$.
Here, $\,\zeta^4-1=(\zeta-1)\zeta^j\,$ for some $\,j$, so 
$\,\zeta^j=\zeta^3+\zeta^2+\zeta+1$.  As before, this gives $\,5\in W_p(m)$.
If $\,p>7$, then $\,d=(p-1)/2\geq 5$, and we have what we want.  Thus we are
only left with the cases $\,p=3,\,5,\,7$.

\ms
If $\,p=3$, we have $\,d=(p-1)/2=1$ (and $\,m=3^k-1\,$ by (4.3)).  In this 
case the desired conclusion is $\,[\,2,\,\infty)_{{\mathbb Z}}\subseteq W_3(m)$, 
which is true since $\,2,\,3\in W_3(m)$.  

\ms
If $\,p=5$, then $\,d=2\in W_5(m)$. In this case we need to show that 
$\,3\in W_5(m)$.  If $\,a_1G=a_2G$, we are done as 
before. Otherwise, one of these cosets must be the identity coset $\,G$
(since $\,[\dot{K}:G]=d=2$), and this implies again that $\,3\in W_5(m)$.

\ms
Finally, we treat the case $\,p=7$. Here we must show that
$\,[\dot{K}:G]=d=3\,$ is in the weight set $\,W_7(m)$.  We first note that:

\bs
\nt (4.4) {\it If $\,(\zeta-1)G=(\zeta^i-1)G=(\zeta^{i+1}-1)G\,$ for some
$\,i\geq 1$, then $\,3\in W_7(m)$. }

\bs
Indeed, if we write $\,\zeta^i-1=(\zeta-1)\zeta^r\,$ and
$\,\zeta^{i+1}-1=(\zeta-1)\zeta^s$, then
$\,\zeta^r=\zeta^{i-1}+\cdots+\zeta+1\,$ and
$\,\zeta^s=\zeta^{i}+\cdots+\zeta+1\,$ imply that
$\,\zeta^s=\zeta^i+\zeta^r$, so $\,3\in W_7(m)$.  Now let $\,C,\,C'\,$ be the 
two nonidentity cosets of $\,G\,$ in $\,\dot{K}$.  By reasonings we have used
before, we may assume that $\,(\zeta^2-1)G=C\,$ and 
$\,(\zeta-1)G=(\zeta^4-1)G=C'$.  Noting that $\,3\,$ is prime to $\,m$
(since $\,\mbox{gcd}(p-1,\,m)=\mbox{gcd}(6,\,m)=2$), we may also assume, in 
view of (4.4), that $\,(\zeta^3-1)G=C$.  Replacing $\,\zeta\,$ by $\,\zeta^3$,
we may further assume that $\,(\zeta^9-1)G=C'$.  Next, note that since
$\,m\geq 6$, it cannot divide $\,10$, so $\,\zeta^{10}\neq 1$.  Thus,
in view of (4.4), we may assume that
$\,(\zeta^5-1)G=C$, and hence that $\,(\zeta^{10}-1)G=C'$.  Now we have
$\,C'=(\zeta-1)G=(\zeta^9-1)G=(\zeta^{10}-1)G$, so $\,3\in W_7(m)\,$ once more
by (4.4).   \qed

\section{The Case $\,\mbox{gcd}(p-1,\,m)=1$}

{\it Throughout this section, we shall assume that} $\,\mbox{gcd}(p-1,\,m)=1\,$
(and as before, $\,m\geq 3$). The analysis of the weight set 
$\,W_p(m)\,$ in this case turns out to require the hardest work.  
  
\bs
The assumption that $\,\mbox{gcd}(p-1,\,m)=1\,$ means that the only 
$\,m\,$th root of unity in $\,{\mathbb F}_p\,$ is $\,1$. Therefore, upon
factoring the polynomial $\,X^m-1\,$ modulo $\,p$, we have
$$  X^m-1=(X-1)g_1(X)g_2(X)\cdots \,,  \leqno (5.1)  $$
where the $\,g_i\,$'s are irreducible monic polynomials in $\,{\mathbb F}_p[X]$, 
each of degree $\geq 2$.  Let $\,\ell:=\mbox{min}\{\mbox{deg}(g_i)\}$.  This 
integer $\,\ell\,$ will play an important role in finding the estimates on 
$\,W_p(m)\,$ in this section, so let us first note a few other 
characterizations of it.

\bs
Recall that the cyclotomic polynomial $\,\Phi_n(X)\in {\mathbb Z}[X]\,$
factors modulo $\,p\,$ into a product of irreducible factors each of
degree given by the order of the element $\,p\,$ in the unit group
$\,U({\mathbb Z}/n\,{\mathbb Z})$ (see, e.g. [Gu]).
Since $\,X^m-1=\prod_{n|m}\Phi_n(X)$, it follows that  $\,\ell\,$ is 
the minimum of the orders of $\,p\,$ in 
$\,U({\mathbb Z}/n\,{\mathbb Z})\,$ for $\,n\,$ ranging over the divisors
of $\,m\,$ greater than 1.  From this, we see that $\,\ell\,$ is also the 
minimum of the orders of $\,p\,$ in $\,U({\mathbb Z}/q\,{\mathbb Z})\,$ for $\,q\,$
ranging over the prime divisors of $\,m$. It is now an easy exercise to
check the following:
$$ \ell = \mbox{min}\,\{e\geq 1:\,\; \mbox{gcd}(p^e-1,\,m)>1\}.  
\leqno (5.2)$$
This simply means that $\,{\mathbb F}_{p^{\ell}}\,$ is the field with the
smallest extension degree over $\,{\mathbb F}_p\,$ which contains an $\,m\,$th 
root of unity other than 1. This can also be verified directly from 
the definition of $\,\ell$.

\bs
For the rest of this section, let  $\,L:={\mathbb F}_{p^{\ell}}$,
$\,m':=\mbox{gcd}(p^{\ell}-1,\,m)$, and let $\,H\,$ be
group of $\,m'\,$th roots of unity in $\,L$.  By (5.2), $\,|H|=m'\geq 2$. 
It will be important to work with the set $\,T:=\mbox{tr}(H)\,$
where ``tr'' denotes the field trace from $\,L\,$ to $\,{\mathbb F}_p$.
The next theorem gives a description of $\,W_p(m)\,$ in terms of 
$\,\ell\,$ and the cardinality $\,t:=|T|\,$ (under the standing assumption 
that $\,\mbox{gcd}(p-1,\,m)=1$).

\bs
\nt {\bf Theorem 5.3.} {\it  Let $\,\ell\,$ and $\,t\,$ be as defined above, 
and let $\,n:=\lceil \frac{p-1}{t-1} \rceil$.  Then
$\,[\,\ell\,n,\,\infty)_{{\mathbb Z}}\subseteq W_p(m')\subseteq W_p(m)$.}

\bs
\nt {\bf Proof.}  Applying the Cauchy-Davenport Theorem
to the subset $\,T\,$ in $\,{\mathbb F}_p$, we have
$\,|2\ast T|\geq \mbox{min}\,\{p,\,2t-1\}$, and inductively
$\,|i\ast T|\geq \mbox{min}\,\{p,\,it-(i-1)\}$.  By the definition of $\,n$,
we have $\,n(t-1)\geq p-1$, so $\,nt-(n-1)\geq p$.  Therefore, 
$\,|n\ast T|=p$.  In particular, for every $\,j\geq 0$, there exists an 
equation $\,t_1+\cdots+t_{n}=-j\in {\mathbb F}_p$, where all $\,t_i\in T$.
Now each $\,t_i\in \mbox{tr}(H)\,$ is a sum of $\,\ell\,$ 
elements of $\,H$, so $\,t_1+\cdots+t_{n}+j\cdot 1=0\,$ \linebreak
is a vanishing sum of $\,m'\,$th roots of unity of weight $\,\ell\,n+j$.  
This shows that
$\,[\,\ell\,n,\,\infty)_{{\mathbb Z}}\subseteq W_p(m')\subseteq W_p(m)$, 
as desired. \qed

\bs
Note that the above theorem is meaningful only if we know that the trace set
$\,T\subseteq {\mathbb F}_p\,$ has at least two elements.  Fortuitously, this
is always the case, according to the following result.

\bs
\nt {\bf Trace Lemma 5.4.} {\it In the notations of} $\,(5.3)$, 
$\,t=|T|\geq 2$.

\bs
The proof of this lemma will be postponed to the last section (\S6).
We shall first assume this lemma and try to get to the main conclusions
of this section.  Note that the larger the trace set $\,T\,$ is, the
better bound on $\,W_{p}(m)\,$ is given by (5.3).  Since $\,t\geq 2\,$ by
(5.4), we have in any case:

\bs 
\nt {\bf Corollary 5.5.} 
$\,[\,\ell\,(p-1),\,\infty)_{{\mathbb Z}}\subseteq W_p(m')\subseteq W_p(m)$.

\bs
Now consider any field $\,K={\mathbb F}_{p^{k}}\,$ containing the group
$\,G\,$ of all $\,m\,$th roots of unity.  Clearly, $\,L\subseteq K$,
and $\,H=\dot{L}\cap G$.  Let $\,d=(p^k-1)/m\,$ and $\,d'=(p^{\ell}-1)/m'$.
We now proceed to the proof of the following, which is a stronger version of
Theorem 1.3 in the case treated in this section.

\bs
\nt {\bf Theorem 5.6.} {\it Assume that} $\,\mbox{gcd}(p-1,\,m)=1$. {\it Then}
$\,[\,d\,',\,\infty)_{{\mathbb Z}}\subseteq W_p(m')\,$ {\it and}
$\,[\,d,\,\infty)_{{\mathbb Z}}\subseteq W_p(m)\,$ {\it except in the following
two special cases:} (A)  $\,d\,'=p-1$; (B) $\,p=2,\;d\,'=3,\;and\;\,m'=5.$
{\it In these special cases, we have} 
$\,[\,d\,'+1,\,\infty)_{{\mathbb Z}}\subseteq W_p(m')\,$ {\it and}
$\,[\,d+1,\,\infty)_{{\mathbb Z}}\subseteq W_p(m)$.

\bs
\nt {\bf Proof.}  Since $\,d\,'=[\dot{L}:H]=[\dot{L}G:G]\,$ divides
$\,d=[\dot{K}:G]\,$ and $\,m'\,|\,m$, it suffices to prove the theorem
for $\,W_p(m')$. Let 
$$\,s=(p^{\ell-1}+\cdots+p+1)/m', \leqno (5.7) $$ 
so that $\,d\,'=s(p-1)$.  {\it First let us treat the special case} (A),
{\it where we have $\,s=1$.}  Here we are supposed to prove that 
$\,[p,\,\infty)_{{\mathbb Z}}\subseteq W_p(m')$.  Since now 
$\,[\dot{L}:H]=p-1\,$ (and $\,H\cap {\mathbb F}_p=\{1\}$), we have
$\,\dot{L}=H\cdot \dot{{\mathbb F}}_p$.  Therefore, fixing a primitive
$\,m'\,$-th root of unity $\,\alpha\in H$, we have
 $\,\alpha-1=b^{-1}\alpha^i\,$ for some integer $\,i\,$ and some $\,b\in 
\dot{\mathbb F}_p$.  For convenience, let us think of $\,b\,$ as an integer
in $\,[1,\,p-1]$.  Multiplying $\,b\cdot 1+(p-b)\cdot 1=0\in {\mathbb F}_p\,$
by $\,\alpha\,$ and using the relation $\,b\alpha=b+\alpha^i$, we get
$$ 0=b\alpha+(p-b)\alpha=b\cdot 1+\alpha^i+(p-b)\alpha,  $$
which is a vanishing sum of weight $\,p+1$.  Multiplying this by $\,\alpha\,$
again and repeating the argument, we get vanishing sums (of $\,m'\,$-th roots
of unity) of weight $\,p+i\,$ for any $\,i>0$.  Coupled with
$\,p\in W_p(m')$, this gives
$\,[\,p,\,\infty)_{{\mathbb Z}}\subseteq W_p(m')$, as desired.  {\it For the rest of
the proof, we may assume that $\,s>1$.} We claim the following:

\bs
\nt {\bf Lemma 5.8.} $\,s>1\,$ {\it implies that\/} $s\geq \ell$,
{\it except perhaps when $\,p=2\,$ and $\,\ell=4,6,8,9$.}

\bs
Thus, leaving aside the four special cases, we have 
$\,d\,'=s(p-1)\geq \ell\,(p-1)$, so the desired conclusion for $\,W_p(m')\,$ 
in (5.6) follows from (5.5).  The four special cases will have to 
be treated later.

\bs
\nt {\bf Proof of (5.8).} We go into the following two cases.

\ms
\nt {\it Case 1.} $\,\ell\,$ {\it is prime}.  We claim that $\,\ell\leq q\,$
for any prime $\,q\,|\,s\,$ (and therefore $\,\ell\leq s$).  In fact, from
(5.7), we get $\,p^{\ell-1}+\cdots+p+1\equiv 0\;(\mbox{mod}\;q)$, so 
$\,p^{\ell}\equiv 1\;(\mbox{mod}\;q)$.  If $\,p\not\equiv 1\;(\mbox{mod}\;q)$, 
then $\,p\,$ has order $\,\ell\,$ in $\,U({\mathbb Z}/q{\mathbb Z})\,$  
(since $\,\ell\,$ is prime), and so $\,\ell\,|\,(q-1)$.  In this case 
$\,\ell\leq q-1<q$.  If $\,p\equiv 1\;(\mbox{mod}\;q)$, then from
$\,p^{\ell-1}+\cdots+p+1\equiv \ell\equiv 0\;(\mbox{mod}\;q)$, we have in
fact $\,\ell=q$.

\ms
\nt {\it Case 2.} $\,\ell\,$ {\it is composite.}  Let $\,q\,$ be the smallest
prime divisor of $\,\ell\,$ and write $\,\ell=qt$.  Then
$\,1<t<\ell\,$ and (5.2) implies that $\,\mbox{gcd}(p^t-1,\,m)=1$.
Since $\,(p^t-1)\,|\,(p^{\ell}-1)$,
we see that $\,p^t-1\,$ divides $\,(p^{\ell}-1)/m'=s(p-1)$.  Thus,
$\,s\geq (p^t-1)/(p-1)$.  We shall now exploit the following elementary
fact which is easy to prove using calculus:

\bs
\nt {\bf Lemma 5.9.} $\,p^x\geq (p-1)x^2+1\,$ {\it for every}
$\,x\in [\,2,\,\infty)_{{\mathbb Z}}\,$ {\it with the exception of $\,p=2\,$ and
$\,x=2,\,3,\,4.$}

\bs
Applying this lemma to $\,x=t$, we get the desired conclusion
$$ s\geq \frac{p^t-1}{p-1}\geq t^2\geq qt=\ell, $$
except when $\,p=2\,$ and $\,t=2,3,4$.  If $\,t=2$, we have $\,q=2\,$ so
$\,\ell=4$.  If $\,t=3$, we have $\,q=2,3$, so $\,\ell=6\,$ or $\,9$.
Finally, if $\,t=4$, we have $\,q=2\,$ so $\,\ell=8$.  This proves (5.8),
but we still have to complete the proof of (5.6) in the four special cases
noted.

\bs
In these cases, $\,d\,'=(2^{\ell}-1)/m'$, so both $\,m',\,d\,'\,$ are odd (and
$>1$).  {\it We may assume that $\,d\,'<m'$}.  (For, if $\,d\,'\geq m'$, we 
have $\,[\,m',\,\infty)_{{\mathbb Z}}\subseteq W_2(m')\,$ since $\,m'\,$ is odd, 
and hence $\,[\,d\,',\,\infty)_{{\mathbb Z}}\subseteq W_2(m')$.)  We simply have to 
check the four outstanding cases individually.

\ms
\nt (1) $\ell=4$. Here $\,2^{\ell}-1=15$, so $\,d\,'=3,\;m'=5$, and we are in
the case (B) of (5.6).  The desired conclusion in this case is 
$\,[\,4,\,\infty)_{{\mathbb Z}}\subseteq W_2(m')$, which is true 
since $\,5=m'\in W_2(m')$.  In fact, by (2.6), we have
$\,W_2(5)=2\,{\mathbb N}+5\,{\mathbb N}=\{0,2\}\cup [\,4,\,\infty)_{{\mathbb Z}}$.  
In particular, $\,3\notin W_2(5)$, so this case is truly exceptional.

\sk
\nt (2) $\ell=6$.  Here $\,2^{\ell}-1=63$, so we have either $\,d\,'=7,\;
m'=9\,$ or $\,d\,'=3,\,m'=21$.  In both cases, 
$\,[\,2,\,\infty)_{{\mathbb Z}}\subseteq W_2(m')\,$ (since $\,2,3\in W_2(m')$),
so there is no problem.  (Actually, in the case $\,d\,'=7$, we are in the 
good case $\,s=d\,'\geq \ell\,$ already.) 

\sk
\nt (3) $\ell=8$. Here $\,2^{\ell}-1=255$, so we have either $\,d\,'=5,\;
m'=51\,$ or $\,d\,'=15,\,m'=17$.  The latter case presents no problems, since 
we are once more in the good case  $\,s=d\,'\geq \ell$.   In the former case,
$\,3\,|\,m'\,$ implies that $\,W_2(m')=[\,2,\,\infty)_{{\mathbb Z}}$, 
so again there is no problem.

\sk
\nt (4) $\ell=9$.  Here $\,2^{\ell}-1=511$, so  $\,d\,'=7,\;m'=73$.  We have
shown in (2.5) that $\,W_2(73)=\{0\}\cup [\,2,\,\infty)_{{\mathbb Z}}$, so
there is no problem.

\bs
This finally completes the proof of Theorem 5.6.

\section{Traces of $\,m\,$th Roots of Unity}

In \S5, we stated without proof the Trace Lemma 5.4, which was crucial for
the proofs of (5.5) and (5.6).  In this section, we return to the trace set 
$\,T=\mbox{tr}(H)$, and offer a general analysis of $\,T\,$  which we believe
to be of independent interest.  The proof of the Trace Lemma is an easy 
by-product of this general analysis.

\bs
The notations (and hypotheses) introduced  at the beginning of \S5 will 
remain in force.  In particular, $\,H\,$ is the group of $\,m'\,$th roots 
of unity in $\,L={\mathbb F}_{p^{\ell}}$, and ``tr'' is the field trace from
$\,L\,$ to $\,{\mathbb F}_p$.  To enumerate the elements in $\,T$, let
$$\,X^{m'}-1=(X-1)h_1(X)\cdots h_r(X)  \leqno (6.1) $$
be the factorization of $\,X^{m'}-1\,$ into (monic) irreducibles over 
$\,{\mathbb F}_p$. Then $\,\mbox{deg}\;h_i\geq \ell\,$ by the definition of 
$\,\ell\,$ (and the fact that $\,(X^{m'}-1)\,|\,(X^m-1)$).  On the other hand,
since $\,L\,$ contains all $\,m'\,$th roots of unity, each $\,h_i(X)\,$ splits 
completely in $\,L$, so $\,\mbox{deg}\;h_i  \leq [L:{\mathbb F}_p]=\ell$. 
Therefore, $\,\mbox{deg}\;h_i=\ell\,$ for all $\,i$. Let
$$ h_i(X)=X^{\ell}-a_iX^{\ell-1}+\,\cdots \,,  \leqno (6.2)  $$
 and let $\,\{\alpha_{ij}\}\,$
be all the roots of $\,h_i(X)\,$ in $\,L$.  For each $\,h_i$, we can identify 
the field $\,{\mathbb F}_p[X]/(h_i(X))\,$ with $\,L\,$ by the correspondence 
$\,\overline{X}\leftrightarrow \alpha_{ij}\,$ (for any $\,j$).
Therefore, $\,\mbox{tr}(\alpha_{ij})=\sum_{k}\,\alpha_{ik}=a_i\,$ for all
$\,i,\;j$.  We have thus
$$ T=\{\mbox{tr}(1),\,a_1,\,\cdots\,,\,a_r\}=\{\ell,\,a_1,\,\cdots\,,\,a_r\},
           \leqno (6.3) $$
with possible duplications. 

\bs
It is now easy to prove the Trace Lemma 5.4, which asserted that 
$\,|T|\geq 2$. Assume, for the moment, that $\,T\,$ is a singleton.  
Then, by (6.3), $\,a_{i}=\ell\,$ for all $\,i$.  Summing all roots 
of the polynomial in (6.1) (and recalling that $\,m'\geq 2$), we get
$$ 0=1+a_1+\cdots+a_r=1+r\,\ell=m'\in {\mathbb F}_p \;, $$
contradicting the fact that $\,m'\,$ is prime to $\,p$.

\bs
The equation (6.3) gives an upper bound $\,|T|\leq 1+r$, and this becomes
an equality iff the elements listed in (6.3) are distinct. {\it  This is the
case, for instance, if} $\,\ell=2$.  To see this, note that the constant
term of each $\,h_i(X)\,$ in (6.1) is $\,1$, since it is an $\,m\,$th 
root of unity in $\,{\mathbb F}_p$, and we are assuming that
$\,\mbox{gcd}(p-1,\,m)=1$.  Therefore, if $\,\ell=2$, we have
$\,h_i(X)=X^2-a_iX+1$.  Since the $\,h_i\,$'s are distinct, so are the
$\,a_i\,$'s, and of course $\,a_i\neq 2\,$ (since otherwise 
$\,h_i(X)=(X-1)^2$).  Therefore, $\,|T|=1+r$, and (5.3) gives the pretty good 
estimate $\,[\,2n,\,\infty)_{{\mathbb Z}}\subseteq W_p(m')\subseteq W_p(m)\,$ 
with $\,n=\lceil \frac{p-1}{r}\rceil$.  For a simple example of this, 
let $\,p=5$, and $\,m=3$.  Here $\,m'=3$, $\,\ell=2$, $\,r=1$, 
$\,L={\mathbb F}_{25}$, and $\,|T|=2$.  By (2.6), $\,W_5(3)\,$ is the set
$$ 3\,{\mathbb N}+5\,{\mathbb N}=\{0,\,3,\,5,\,6\}\cup [\,8,\,\infty)_{{\mathbb Z}}. $$
Since $\,2n=8$, the conclusion in (5.3) is sharp here.

\bs
In general, $\,t:=|T|\,$ may be less than $\,1+r$, since there may be 
duplications among the elements of $\,T\,$ listed in (6.3).  For an example
where interesting duplications occur, take $\,p=7\,$ and $\,m=19$.
Here $\,m'=19\,$ and $\,\ell=3,\;\,r=6$. {\it Mathematica\/} gives a 
factorization 
\begin{eqnarray*}
 X^{19}-1&=&(X-1)(X^3+2X+6)(X^3+4X^2+X+6)(X^3+4X^2+4X+6)  \\
         & &  \mbox{} \cdot(X^3+5X^2+6)(X^3+3X^2+3X+6)(X^3+6X^2+3X+6)
              \in {\mathbb F}_{7}[X]. 
\end{eqnarray*}
Since $\,-4=3=\mbox{tr}(1)\,$ in $\,{\mathbb F}_{7}$, $\,T\,$ has only five 
(two less than $\,1+r=7$) distinct elements $\,\{0,\,1,\,2,\,3,\,4\}$.   
In this case, the number $\,n\,$ in (5.3) is 
$\,\lceil \frac{7-1}{5-1} \rceil =2$, and 
(5.3) shows that $\,[\,6,\,\infty)_{{\mathbb Z}}\subseteq W_{7}(19)$.  Note that,
in spite of the trace duplications, this is still much sharper than what is 
given in (5.6).

\bs
In general, we cannot hope to improve upon the lower bound $\,|T|\geq 2$.  
For one thing, $\,{\mathbb F}_p\,$ may have only two elements to begin with.
Also, we may have $\,r=1$, in which case (5.4) and (6.3) show that $\,|T|=2$.
Even if $\,p\geq 3\,$ and $\,r\geq 2$, there are many cases in which $\,T\,$
is just a doubleton.  Let us illustrate the situation $\,r=2\,$ by taking
$\,m\,$ to be an odd prime $\,q\,$ (so that $\,m'=q\,$ too), and assuming 
that $\,p\,$ is also odd and has order $\,\ell=(q-1)/2\,$ in the group 
$\,U({\mathbb Z}/q\,{\mathbb Z})$.  In this case, $\,r=2$, and (6.1) becomes
$$\,X^{q}-1=(X-1)h_1(X)h_2(X)  \leqno (6.4) $$
where $\,h_1,\;h_2\,$ are monic irreducible (over $\,{\mathbb F}_p\,$) of degree
$\,\ell$.  Following a standard notation in number theory, let us define
$\,q^*\,$ to be $\,q\,$ if $\,q \equiv 1\,(\mbox{mod}\;4)$, and
$\,q^*\,$ to be $\,-q\,$ if $\,q \equiv 3\,(\mbox{mod}\;4)$.  Then the
size of the trace set $\,T=\mbox{tr}(H)\,$ is determined as follows.

\bs
\nt {\bf Proposition 6.5.} {\it Under the above assumptions $\,|T|=2\,$
iff $\,p\,|\,(q^*-1)\,$ (and $\,|T|=3\,$ otherwise).}

\bs
\nt {\bf Proof.}  Let $\,E={\mathbb Q}(\zeta)\,$ where $\,\zeta =e^{2\pi i/q}$,
and fix a generator $\,t\,$ of $\,U({\mathbb Z}/q\,{\mathbb Z})$. Then 
$\,\sigma:\,\zeta \mapsto \zeta^t\,$ is a generator for 
$\,\mbox{Gal}(E/{\mathbb Q})$, and
$\,\sigma^2:\,\zeta \mapsto \zeta^{t^2}\,$ is a generator for 
$\,\mbox{Gal}(E/F)$, where $\,F\,$ is the fixed field $\,E^{\sigma^2}$.  Note 
that $\,X^{q-1}+\cdots+X+1\,$ factors into $\,f(X)g(X)\,$ over $\,F[X]$, where 
$$ f(X)=X^{\ell}-a\,X^{\ell-1}+\cdots  \;\;\; \mbox{and} \;\;\;
   g(X)=X^{\ell}-b\,X^{\ell-1}+\cdots\,  $$
are, respectively, the minimal polynomials of $\,\zeta\,$ and $\,\zeta^t\,$
over $\,F$.  We have
$$ a-b=\mbox{tr}_{E/F}(\zeta)-\mbox{tr}_{E/F}(\zeta^t)
      =\sum_{j=0}^{\ell-1}\,\zeta^{t^{2j}}-\sum_{j=0}^{\ell-1}\,\zeta^{t^{2j+1}} \in F.  $$
This is precisely the quadratic Gauss sum (with respect to the Legendre
character on $\,{\mathbb F}_q\,$), so by [IR:~(8.2.2)], $\,a-b=\sqrt{q^*}\,$.  
(Gauss showed that the $\,\sqrt{q^*}\,$ here is the one taken in the upper 
half plane if $\,q\equiv 3\,(\mbox{mod}\;4)$, but this will not be needed in 
the following.)  Since we also have $\,a+b=-1$, it follows that 
$$  a=(\sqrt{q^*}-1)/2\,, \;\;\;\;  b=-(\sqrt{q^*}+1)/2.    $$
Incidentally, this proves the well-known fact that 
$\,F={\mathbb Q}(\sqrt{q^*}\,)$.

\mk
Let $\,R\,$ be the ring of algebraic integers in $\,F$.  Since $\,p\,$
is unramified in $\,E$, it is also unramified in $\,F$, so 
$\,pR=\mathfrak{p}\,\mathfrak{p}'\,$, where $\,\mathfrak{p},\;\mathfrak{p}'\,$ are distinct
prime ideals of $\,R$, both of residue degree 1. Identifying
$\,R/\mathfrak{p}\,$ with $\,{\mathbb F}_{\it p}$, we may take the polynomials 
$\,h_1,\;h_2\,$ in (6.4) to be $\,\overline{f}\,$ and $\,\overline{g}$, 
where ``bar'' means reduction modulo $\,\mathfrak{p}$.  In particular,
$\,T=\{\ell,\,\overline{a},\,\overline{b}\}\,$ by (6.3).  Here the two elements
$\,\overline{a},\;\overline{b}\,$ are {\it always\/} different (for otherwise
$\,\mathfrak{p}\,$ would contain $\,(a-b)^2=q^*\,$ as well as $\,p$).
Therefore, $\,T\,$ will have only two elements iff $\,\mathfrak{p}\,$ also
contains 
$$ 4(\ell-a)(\ell-b)=(q-\sqrt{q^*}\,)\,(q+\sqrt{q^*}\,)=q^2-q^*=q^*(q^*-1).$$
Since $\,q^*\notin \mathfrak{p}$, this happens iff $\,q^*-1\in \mathfrak{p}$, that is,
iff $\,p\,|\,(q^*-1)$, as claimed. \qed

\bs
\nt {\bf Corollary 6.6.} {\it Let $\,q=2\ell+1\,$ and $\,p=2\ell\,'+1\,$
be distinct primes such that the order of $\,p\,$ is $\,\ell\,$
modulo $\,q$.  If $\,p\,|\,(q^*-1)$, then
$\,[\,2\ell\ell\,',\,\infty)_{{\mathbb Z}}\subseteq W_p(q)$.  Otherwise,
$\,[\,\ell\ell\,',\,\infty)_{{\mathbb Z}}\subseteq W_p(q)$.}

\bs
\nt {\bf Proof.} This follows from (5.3) and (6.5), since the number $\,n\,$
in (5.3) is $\,2\ell\,'\,$ in the first case, and $\,\ell\,'\,$ in the
second case. \qed

\bs
\nt {\bf Example 6.7.} Let $\,q=11\,$ (with $\,q^*=-11$).  Then the primes 
$\,3\,$ and $\,5\,$ both have order $\,\ell=(q-1)/2=5\,$ modulo $\,q$, and
according to {\it Mathematica\/}:
$$X^{11}-1=(X-1)(X^5+X^4+2X^3+X^2+2)(X^5+2X^3+X^2+2X+2)\in {\mathbb F}_{3}[X],$$
\vspace{-.6cm}
$$X^{11}-1=(X-1)(X^5+2X^4+4X^3+X^2+X+4)(X^5+4X^4+4X^3+X^2+3X+4)\in {\mathbb F}_{5}[X].   $$

\vspace{.3cm}
\nt Thus, for $\,p=3$, $\,T=\{\ell,\,-1,\,0\}=\{2,\,0\}\,$ in 
$\,{\mathbb F}_3$.  This is consistent with (6.5) since $\,p=3\,$ divides 
$\,q^*-1=-12$.  Here, $\,\ell=5,\;\ell\,'=1$, so (6.6) gives 
$\,[\,10,\,\infty)_{{\mathbb Z}}\subseteq W_3(11)$. (In fact, from $\,0,\,2\in T$,
we see easily that $\,5,\,6\in W_3(11)$, and so $\,8,\,9\in W_3(11)\,$ also.)
On the other hand, if we choose $\,p=5$, then
$\,T=\{\ell,\,-2,\,-4\}=\{0,\,3,\,1\}\,$ in $\,{\mathbb F}_{5}$, 
consistently with (6.5) since $\,p=5\,$ does not divide $\,q^*-1=-12$.  
Here, $\,\ell=5,\;\ell\,'=2$, so (6.6) gives again
$\,[\,10,\,\infty)_{{\mathbb Z}}\subseteq W_5(11)$, and $\,0,\,1,\,3\in T\,$
show further that $\,5,\,7,\,9\in W_{5}(11)$.

\bs
The arguments in the proof of (6.5) can be generalized.  However, if the
order of $\,p\,$ modulo $\,q\,$ is smaller than $\,(q-1)/2\,$ (in other
words $\,r>2$), the computations of the trace elements in $\,T\,$ will 
involve Gaussian sums with (higher) character values as coefficients.  
We shall not go into this analysis here.  We should point out, however, 
that if $\,q\,$ is {\it fixed\/}, then the prime ideal method 
(in characteristic 0) used in the proof of (6.5) will suffice to show that 
the upper bound $\,|T|\leq 1+r\,$  becomes an {\it equality\/} for sufficiently
large $\,p$.  Therefore, by (5.3), 
$\,[\,\ell\,n,\,\infty)_{{\mathbb Z}}\subseteq W_p(q)\,$ with
$\,n=\lceil \frac{p-1}{r} \rceil\,$, for sufficiently large $\,p$.

\bs

\bs

\bs

\end{document}